\documentclass[12pt]{article}

          \usepackage{graphics}

\usepackage{latexsym,amsmath,amssymb,amsthm}
\newtheorem{tw}{Theorem}[section]

\newtheorem{stw}[tw]{Proposition}
\newtheorem{wn}[tw]{Corollary}

\theoremstyle{definition}
\newtheorem{df}[tw]{Definition}
\newtheorem{prz}[tw]{Example}

\theoremstyle{remark}
\newtheorem{uw}[tw]{Remark}

\numberwithin{equation}{section}

\newcommand{\zaw}{\preccurlyeq}
\newcommand{\row}{\thickapprox}

\newcommand{\ba}{\begin{array}}
\newcommand{\ea}{\end{array}}
\newcommand{\beq}{\begin{equation}}
\newcommand{\eeq}{\end{equation}}
\newcommand{\C}{\mathbb{C}}
\newcommand{\cc}{\mathcal{C}}
\newcommand{\cZ}{\mathcal{Z}}

\newcommand{\f}{\varphi}
\newcommand{\ov}{\overline}

\newcommand{\al}{\alpha}

\newcommand{\N}{\mathbb{N}}
\newcommand{\Oh}{\mathcal{O}}

\newcommand{\R}{\mathbb{R}}

\newcommand{\ti}{\tilde}

\newcommand{\sk}{\fin(T)}
\newcommand{\ski}{\fin(I)}

\newcommand{\ani}{\al \in \fin(I)}
\newcommand{\On}{\mathcal{O}_{x_0}(K^T)}
\newcommand{\Op}{\mathcal{O} }
\newcommand{\A}{\mathcal{A} }
\newcommand{\B}{\mathcal{B} }
\newcommand{\Os}{\mathcal{O}_{x_{0S}}(K^S)}

\newcommand{\zal}{Z(\al )}
\newcommand{\pp}{\left(\Pi^T_S\right)^{-1}}

\newcommand{\cf}{\mathcal{F}_K^T}

\newcommand{\rp}{\{x_0\}_{\sk}}

 \DeclareMathOperator{\fin}{Fin}

    \title{Systems of germs and
        theorems of zeros in  infinite-dimensional spaces}
   \author{Dorota Mozyrska \& Zbigniew Bartosiewicz}

          \date{}

\begin{document}

\maketitle
        \begin{abstract}

            Systems of germs of sets in infinite-dimensional spaces are introduced
            and studied. Such a system corresponds to a local zero-set of an
            ideal of the ring of analytic functions of infinite number of
            variables. Conversely, this system of germs defines the ideal of
            germs of analytic functions vanishing on it. A theorem of zeros is
            proved, stating that this ideal is the radical (in the complex
            case) or real radical (in the real case) of the initial
            ideal.\\
        \end{abstract}
          \maketitle

          \section{Introduction}\label{intro}

          In 1952 S.~Lang \cite{La} extended Hilbert's Nullstellensatz to
polynomials of infitely many variables.  On the other hand,
W.~R\"{u}ckert \cite{Ruc,Ru}, in 1932, instead of polynomials took
germs of complex analytic functions of finitely many variables;
his Nullstellensatz involved germs of complex analytic sets. The
real case for polynomials of finitely many variables was
independently solved by J.-L.~Krivine \cite{Kr}, D.W.~Dubois
\cite{Du} and J.-J.~Risler \cite{Ri0} in the sixties (of the last
century). In their theorem of zeros, the ordinary radical of an
ideal, used in the complex version, had to be changed for the real
radical. Finally, in 1976 Risler \cite{Ri,Ru} proved
finite-dimensional real analytic counterpart of Hilbert's
Nullstellensatz. (See \cite{ABR} for a more abstract real theorem
of zeros.) It is important to notice that the ring of the germs of
analytic functions of finitely many complex or real variables (at
some point) is Noetherian. Hence, if $I$ is an ideal of this ring,
then the germ of zero-set of $I$ is well defined as $I$ is
finitely generated (see e.g. \cite{GR}).

In this paper we study the infinite-dimensional analytic (complex
and real) case, where an analytic function depends (like Lang's
polynomials) only on a finite number of variables (is
\emph{finitely presented}). But, as the number of all variables is
infinite, the ring of the germs of such functions is no longer
Noetherian and the germs of the zero sets of ideals cannot be
defined in the standard way. Moreover, as we have shown in
\cite{MB3}, there is no topology in the infinite-dimensional space
of all complex or real sequences that would give required
properties of the germs of sets. We have been interested there in
``local'' solutions of infinitely many analytic equations in
infinitely many variables. Such equations describe, for instance,
indistinguishable states of a (control) system with output and are
related to \emph{observability} and \emph{local observability} of
the system (see e.g. \cite{Ba1,Ba2} for the finite-dimensional
case and \cite{Mo,MB3,MB4} for the infinite-dimensional one). In
particular, it is important in local observability whether such
equations have locally only one solution (which is the point at
which we localize the system and the equations).

Instead of using topology to define the germ of a set, we consider
special families of finite dimensional set-germs (\emph{systems of
germs}) which approximate in some sense what we want to be an
infinite-dimensional set-germ. Systems of germs give rise to a
concept of multigerm --- the equivalence class of such systems
under a natural equivalence relation. We show that it is the right
language in this infinite-dimensional world. We can manipulate
with multigerms in the same way as we do with finite-dimensional
germs. In particular, we can define multigerm of zeros
corresponding to an ideal of the ring of germs of finitely
presented analytic functions and, conversely, zero ideal of a
multigerm. We consider the real and the complex cases. The main
result of this paper consists of real and complex theorems of
zeros, where we show that the real or ordinary radical of an ideal
consists exactly of the germs of finitely presented analytic
functions (real or complex, respectively) that vanish on the
multigerm of zeros (real or complex again) of the ideal.

We omit here many proofs which are either straightforward or
similar to proofs of earlier statements. Instead, we provide
several examples which give the flavor of the theory. They concern
mostly real functions as the real case is more interesting and
closer to applications. Some of them can be found in D.~Mozyrska's
Ph.D. thesis \cite{Mo}, where local observability of
infinite-dimensional dynamical systems is studied in detail.

          \section{Preliminaries and notation}\label{Prelim}
Let $M$ be a topological space, $x_0\in M$ and $A\subseteq M.$
Then by $A_{x_0}$ we shall denote the germ of the set $A$ (or a
set-germ) at $x_0$.

Whereas the union and the intersection of finite number of
set-germs  are well-defined, these operations on an infinite
collection of set-germs are not necessarily well-defined.

The germ at $x_0$ of the empty set will be called the
\textit{empty set-germ} and the germ of the whole space
--- the \textit{full set-germ} (its representatives are
neighborhoods of $x_0$).  A germ which is not the empty set-germ
will be called a \textit{proper set-germ}.

\begin{df}
Let $M$ and $N$ be topological spaces. Let $h:M\longrightarrow N$
be continuous, $x_0\in M$, $y_0= h(x_0).$ Then by the
\textit{inverse image (at $x_0$) of a germ $A_{y_0}$ with respect
to $h$} we will mean the germ at $x_0$ of the inverse image of a
representative of $A_{y_0},$ i.e.
\begin{equation} h_{x_0}^{-1}(A_{y_0}):=h^{-1}(A)_{x_0}.\label{pobraz}\end{equation}
\end{df}
\begin{stw}
Let  $M,N,P$ be topological spaces,
 $x_0\in M, y_0\in N, z_0\in P$, $A\subset P.$ Let  $g,f$ be continuous mappings,
$g:M\longrightarrow N$,  $f:N\longrightarrow P$,  and $y_0=
g(x_0)$,  $z_0=f(y_0)$.  If $h=f\circ g$, then
$h^{-1}_{x_0}(A_{z_0})=g_{x_0}^{-1}(f^{-1}_{y_0}(A_{z_0})).$

\end{stw}

>From now on let $K=\R$ or $K=\C$.
     If  $x_0\in M$ and $f$ is a $K$-valued function defined
    on some neighborhood of $x_0,$ then by $f_{x_0}$ we shall denote
the germ of the function $f$ (a function-germ) at $x_0$. If $\f$
is a function-germ at $x_0$, then a function $f:M\rightarrow K$
such that $f_{x_0}=\f$ is called a representative of $\f$.

\begin{df}\label{cofk}
Let $M,N$ be topological spaces and $x_0\in M, y_0\in N.$ Let
$g:~M\longrightarrow N$ be a continuous mapping such that
$y_0=g(x_0)$ and $\f$ be the germ of a $K$-valued function at
$y_0.$ Then we define the \textit{pullback of the germ $\f$ with
respect to the map $g$} in the following way:
\begin{equation} g^*_{x_0}(\f):=(f\circ g_{|U})_{x_0},
            \label{cofniecie_kielka} \end{equation}
where $f_{y_0}=\f$ and $f:V\longrightarrow U,$ $V\subset N,$
$U=g^{-1}(V)\subset M.$
\end{df}

In applications we will require that $g^*_{x_0}$ be injective.

\begin{stw}\label{cof2}
Let $g$ be a continuous and open mapping of topological spaces $M$
and $N$. Let $x_0\in M$ and $y_0=g(x_0)$. Then $g^*_{x_0}$ is a
monomorphism from the algebra of germs of all functions at $y_0$
to  the algebra of germs of all functions at $x_0$.
\end{stw}

Let $M=K^n$ and $x_0\in M.$ By $\mathcal{O}_{x_0}(M)$ we denote
the ring of germs of $K$-valued analytic functions at $x_0.$ The
ring $\mathcal{O}_{x_0}(M)$ is Noetherian and local (\cite{Ru}).
Its only maximal ideal will be denoted by $m^n_{x_0}.$ If $I$ is
an ideal of $\Oh_{x_0}(M)$, then $Z(I)$ will denote the zero
set-germ at $x_0$ of $I$. Let $J(Z(I))$ be the ideal in
$\Oh_{x_0}(M)$ of all germs of analytic functions that vanish on
$Z(I).$

Let $P$ be a commutative ring with a unit and $I$ be an ideal of
$P$ (which will be denoted by $I\lhd P$). Then the \textit{real
radical of $I$}, denoted by $\sqrt[\R]{I},$ is the set of all
$a\in P$ for which there are $m\in \N, k\in \N \cup \{0\}$ and
$b_1,\ldots b_k \in P$, such that
\[a^{2m}+b_1^2+\cdots +b_k^2 \in I. \]
Then $I\subseteq \sqrt{I}\subseteq \sqrt[\R]{I}$, where $\sqrt{I}$
denotes the ordinary radical of $I$.
\begin{stw}\label{homom}
Let $f: P\rightarrow R$ be a homomorphism of rings. Then
$f^{-1}(\sqrt{I})=\sqrt{f^{-1}(I)}.$
\end{stw}
In case of real radicals we have a weaker statement:
\begin{stw}\label{homom1}
Let $f: P\rightarrow R$ be a homomorphism of rings. Then
$\sqrt[\R]{f^{-1}(I)}\subseteq f^{-1}(\sqrt[\R]{I}).$
\end{stw}

\begin{uw} Real and complex analytic theorems of zeros may
now be stated as follows,({\rm \cite{Ri,GR,Ru,Ruc}}):
\begin{enumerate}
\item Let $x_0\in \R^n$. If $I$ is an ideal of $\Oh_{x_0}(\R^n)$
then
 $J(Z(I))=\sqrt[\R ]{I}.$
\item Let $z_0\in \C^n$. If $I$ is an ideal of $\Oh_{z_0}(\C^n)$
then
 $J(Z(I))=\sqrt{I}.$
\end{enumerate}
\end{uw}

\section{Germs of finitely presented functions}

Let $T$ be an arbitrary nonempty set and $K=\R$ or $K=\C$.
Consider the product space $\prod\limits_{t\in T} K=K^T=\{x:
T\rightarrow K\}$ with the product topology.
 We denote $x_t:=x(t)$
for $t\in T$ and $x\in K^T.$

>From now on we assume that $T$ is an infinite set. By $\sk$ we
denote the set of all finite nonempty subsets of $T$. Let
$S\in\sk.$ Then $K^S:=\{x_S=x_{|S}:S\rightarrow K, x\in K^T\}.$ By
$\Pi^{S_2}_{S_1}$ we denote the projection
$\Pi^{S_2}_{S_1}:K^{S_2}\rightarrow K^{S_1}: x_{S_2}\mapsto
x_{S_1},$ where $S_1\subseteq S_2\subseteq T.$ Let $x_0\in K^T$.
Then $\Pi^{S_2}_{S_1}(x_{0S_2})=x_{0S_1}$ and
$\Pi^T_S(x_0)=x_{0S}$. The projection $\Pi^{S_2}_{S_1}$ is a
continuous and open mapping of topological spaces $K^{S_2}$ and
$K^{S_1}$.

The sets of the form $V=(\Pi^T_S)^{-1}(U)$, where $S\in\sk$ and
$U\subseteq K^S$ is open, form the basis of the product topology
of $K^T.$
\begin{df} Let $V$ be an open set in $K^T$ with the product topology.
We say that a function $f:V\longrightarrow K$ is \textit{finitely
presented on $K^T$} if there are $S\in\sk$ and a function
$\overline{f}: \Pi^T_S(V)\longrightarrow K$ such that
$f=\overline{f}\circ(\Pi^T_S)_{|V}.$ The function $\ov{f}$ is
called a \textit{representing function of the function $f$}. We
also say that \textit{$f$ depends on a finite number of variables}
indexed by $S$ or that \textit{$S$ is indexing $f$.} We say that
$f$ is \textit{analytic} if $\ov{f}$ is analytic. By
$\mathcal{F}^T_K$ we denote the family of all analytic finitely
presented functions on $K^T$.
\end{df}
Observe that if $S$ is indexing $f$ then any $S'$ such that
$S'\supseteq S$ is also indexing $f$ and any $S\in \fin(T)$ is
indexing any constant function.

\begin{prz}\cite{Ban}\mbox{}\\
Consider the space $\R^{\N}$ with the product topology. A function
$f:\R^{\N}\longrightarrow \R$ is linear and continuous if and only
if there are $k\in\N$ and $a_1,\ldots,a_k\in\R$ such that for each
$x\in\R^{\N}$, $f(x)=\sum\limits_{i=1}^ka_ix_i$. Hence, linear
continuous functionals on $\R^{\N}$ are finitely presented.
\end{prz}

Let $f\in \cf$ and $x_0\in K^T.$ Then we may consider the germ of
$f$ at $x_0$ in the standard way in the topological space $K^T$
with the product topology. The collection of germs at $x_0$ of
functions from $\cf$ forms a commutative ring with a unit (the
germ of the constant function equal 1). It is denoted by $\On$. We
say that $S$ is indexing a germ $\f\in \On$ if $S$ is indexing
some representative of $\f.$

Let $x_0\in K^T$ and $S_1\subset S_2\subseteq T$. Let us consider
the rings $\Op_{x_{0S_1}}(K^{S_1})$, $\Op_{x_{0S_2}}(K^{S_2})$, of
germs of analytic functions, respectively at points $x_{0S_1},
x_{0S_2}$.

Let us consider a monomorphism of the ring
$\Op_{x_{0S_1}}(K^{S_1})$ into the ring $\Op_{x_{0S_2}}(K^{S_2})$:
\begin{equation} M^{S_2}_{S_1}: \Op_{x_{0S_1}}(K^{S_1})\ni \f\longmapsto
 (\Pi^{S_2}_{S_1})^*_{x_{0S_2}}(\f)
                \in \Op_{x_{0S_2}}(K^{S_2}).
\label{monomorfizm}
\end{equation}
It is a particular case of Definition~\ref{cofk}. Using the
monomorphism defined by (\ref{monomorfizm}) we may identify the
ring $\Op_{x_{0S_1}}(K^{S_1})$ with the subring of the ring
$\Op_{x_{0S_2}}(K^{S_2})$ consisting of germs of functions that do
not depend on variables with indices from $S_2\backslash S_1$.

\begin{stw}
Let $x_0\in K^T.$ Then
\[\On=\bigcup\limits_{S\in
\fin(T)} M^T_S\left(\Op_{x_{0S}}(K^{S})\right).\]
\end{stw}
\begin{proof} Immediately we have that $\bigcup\limits_{S\in \fin(T)}
M^T_S\left(\Op_{x_{0S}}(K^{S})\right)\subseteq \On$ because of the
definition of $M^T_S.$ Now let $\f\in \On$ and $f:
V\longrightarrow K$ be a representative of $\f$. Then there is
$S\in\sk$ and $\ov{f}:\Pi^T_S(V)\longrightarrow K$ such that
$f=\ov{f}\circ(\Pi^T_S)_{|V}$ and
$\f=f_{x_0}=\left(\ov{f}\circ(\Pi^T_S)_{|V}\right)_{x_0}=
\left(\Pi^T_S\right)^*(\ov{f}_{x_{0S}}).$ Hence $\f\in
M^T_S\left(\Op_{x_{0S}}(K^{S})\right).$
\end{proof}

\begin{stw}\label{wystf}
A set $S\in \fin(T)$ is indexing $\f\in \On$ if and only if
$M^T_S\left((M^{T}_{S})^{-1}(\f)\right)=\f.$
\end{stw}

\begin{stw}\label{przeciwobraz}
Let $S_1\subseteq S_2\subseteq T$ and
$I\lhd\Op_{x_{0S_2}}(K^{S_2}).$  Then
$\left(M^{S_2}_{S_1}\right)^{-1}(I):=\{\f\in\Op_{x_{0S_1}}(K^{S_1}):
M^{S_2}_{S_1}(\f)\in I\}$ is an ideal in the ring
$\Op_{x_{0S_1}}(K^{S_1}).$
\end{stw}
\begin{df}
Let $I=(\f_1,\ldots,\f_k)\lhd\On$ and let  $S_i$ be indexing
$\f_i.$ Then $S=\bigcup\limits_{i=1}^k S_i$ is called to
\textit{be indexing $I.$}
\end{df}
\begin{stw}\label{kiel}
Let $I\lhd\On$ be finitely generated and let $S$ be indexing $I.$
Then
$\left[M^T_S\left((M^{T}_{S})^{-1}(I)\right)\right]\cdot\On=I.$
\end{stw}

\begin{stw}\label{kielek_odwr}
Let  $\f\in\On$. The following conditions are equivalent:
\begin{enumerate}
\item  $\f$ is invertible in $\On$. \item  $\exists S\in\sk: \
(M^{T}_{S})^{-1}(\f)$ is invertible in $\Os.$ \item $\f(x_0)\neq
0.$
\end{enumerate}
\end{stw}

\begin{wn} \label{idmaks}\mbox{}\\
 Consider the ideal $m_{x_0}=(x_t-x_{0t})_{t\in T}$ of $\On$ and
let $\f \in \On$.
\begin{enumerate}
\item $\f \ \mbox{is not invertible}\ \Leftrightarrow \f \in
m_{x_0}.$ \item  $\On$ is a local ring with the maximal ideal
$m_{x_0}=(x_t-x_{0t})_{t\in T}$.
\end{enumerate}
\end{wn}
By $\ski$ we denote the set of all finite nonempty subsets of
$I\lhd\On$. Let $\al\in \ski,$ where $I\lhd\On.$ Then $(\al)$
denotes the ideal generated by all the elements of $\al$.

\begin{stw}\label{pot}
Let $I$ be an ideal of $\On$. Then
\[\sqrt{I}=m_{x_0}\Leftrightarrow \forall S\in \sk \ \exists \al\in \ski:
 \sqrt{(\al)}=(x_t-x_{0t})_{t\in S}.\]
\end{stw}

Observe that for the real radical  Proposition~\ref{pot} does not
 hold.

\begin{prz}
Let $x_0=0,$ $T=\N$ and the ideal $I$ of $\mathcal{O}_0(\R^{\N})$
be generated by germs $x_1^2+(x_2-x_3)^2, x_3^2+(x_4-x_5)^2,
\ldots , x_{2k+1}^2+(x_{2k+2}-x_{2k+3})^2,\ldots$. Then
$\sqrt[\R]{I}=m_{x_0}$ but for each $\al \in    \fin(I)$ and
$S\in\sk$: $\sqrt[\R]{(\al)}\neq (x_t-x_{0t})_{t\in S}$.
\end{prz}

As the ring $\On$ is not Noetherian, the zero set-germ of an ideal
may not be well defined. In Section~\ref{roz5} we will define the
zero-system of an ideal $I\lhd\On$ using zero set-germs of finite
subsets of $I.$

Let $\al=\{\f_1,\ldots ,\f_k\}\subset \On.$ Then we define the
zero set-germ of $\al$ in a standard way: $\zal=Z(\f_1)\cap\ldots
\cap Z(\f_k),$ where $Z(\f)=Z(\ti{\f})_{x_0}$ for some
representative of $\f.$ Let $S\in \fin(T)$ be indexing $\al.$ Then
$(M^T_S)^{-1}(\al)=\{(M^T_S)^{-1}(\f_1),\ldots,
(M^T_S)^{-1}(\f_k)\}\subseteq \Os$ and
$Z\left((M^{T}_{S})^{-1}(\al)\right)$ is the set-germ in $K^{S}.$

\begin{stw} \label{rzut}
Let $\al=\{\f_1,\ldots,\f_k\}\subset \On$ and let $S\in \sk$ be
indexing $\al.$ Then
$(\Pi^{T}_{S})^{-1}\left(Z\left((M^{T}_{S})^{-1}(\al)\right)\right)=Z(\al).$
\end{stw}
\begin{proof}
Let for each $i=1,\ldots,k:$ $(f_i)_{x_0}=\f_i$ and
$f_i:V_i\rightarrow K.$ Then there are functions
$\ov{f}_i:\Pi^T_{S}(V_i)\rightarrow K$ such that
$f_i=\ov{f}_i\circ (\Pi^T_{S})_{|V_i}.$ Let $y\in K^T.$ Then:
$y\in
(\Pi^{T}_{S})^{-1}\left(Z(\ov{f}_1,\ldots,\ov{f}_k)\right)\Leftrightarrow
\Pi^{T}_{S}(y)\in Z(\ov{f}_1,\ldots,\ov{f}_k)\Leftrightarrow
\forall i=1,\ldots,k: \ov{f}_i(\Pi^{T}_{S}(y))=0 \Leftrightarrow
\forall i=1,\ldots,k: f_i(y)=0\Leftrightarrow y\in
Z(f_1,\ldots,f_k).$ Hence
$(\Pi^{T}_{S})^{-1}\left(Z(\ov{f}_1,\ldots,\ov{f}_k)\right)_{x_0}=Z(f_1,\ldots,f_k)_{x_0}.$
\end{proof}
\begin{wn}\label{3.15a}
Assume that $\al,\beta \in \fin (\On).$ If $\zal\subseteq
Z(\beta)$ then there is $S\in \fin(T)$ such that
$Z\left((M^{T}_{S})^{-1}(\al)\right)\subseteq
Z\left((M^{T}_{S})^{-1}(\beta)\right).$
\end{wn}
\begin{proof}
Let $S$ be indexing $\al\ \cup \ \beta.$ Then by
Proposition~\ref{rzut} we have that:
$(\Pi^{T}_{S})^{-1}\left(Z\left((M^{T}_{S})^{-1}(\al)\right)\right)\subseteq
(\Pi^{T}_{S})^{-1}\left(Z\left((M^{T}_{S})^{-1}(\beta)\right)\right).$
As $\Pi^{T}_{S}$ is surjective we get
$Z\left((M^{T}_{S})^{-1}(\al)\right)\subseteq
Z\left((M^{T}_{S})^{-1}(\beta)\right).$
\end{proof}

\section{Systems of germs}

By a directed set of indeces we mean an ordered pair
$(\Lambda,\ll),$  where $\Lambda$ is an arbitrary set and $\ll$ is
a transitive relation in $\Lambda$ that  satisfies the
Moore-Smith's condition,\cite{Ras}:
\begin{equation}
\forall \alpha\in\Lambda \ \forall \beta \in \Lambda \ \exists
\gamma \in\Lambda : \ (\alpha \ll \gamma \ \wedge \
\beta\ll\gamma). \label{MS}
\end{equation}

\begin{df}\label{definicja_systemu}
Let $x_0\in K^T$ and $(\Lambda,\ll)$  be a directed set of
indeces. Then we define a \textit{system of germs at $x_0$} as a
set of set-germs at $x_0$: $\{A^{\al},\al \in \Lambda\},$ such
that $\forall \al\in\Lambda$ $\exists S\in \fin(T)$ $\exists \
A\subseteq K^S :$ $A^{\al}=\left(\Pi^T_S\right)^{-1}(A_{x_{0S}})$
and
\[
\forall (\al,\beta \in \Lambda): \ \beta \ll \al \Longrightarrow \
                 A^{\al}\subseteq A^{\beta}.\]
\end{df}

\begin{df} \label{rodzinapunktowa}\mbox{}
\begin{enumerate}
\item A system $\{A^{\al},\al \in \Lambda\}$ such that for every
$\al,\beta \in\Lambda, \ A^{\al}=A^{\beta}$ will be called a {\it
constant system.} \item  The system of germs $\{A^S, S\in \sk\}$
such that for each $S$: $A^S=\pp(x_{0S})$ will be called the {\it
system of point-germs at $x_0$} and will be denoted by
$\{x_0\}_{\sk}$.

\item The system $\{A_S, S\in\sk\}$ such that for each $S$: $A_S$
is the germ of the space $K^T,$ will be called the {\it system of
full germs at $x_0$.} (It is a constant system).

\item We say that a system $\{A^{\al},\al\in\Lambda\}$ is
\textit{trivial} if there exists $\al\in\Lambda$ such that
$A^{\al}=\emptyset$. Then $\alpha\ll\beta \Rightarrow
A^{\beta}=\emptyset$.

\item The system of germs at $x_0$: $\{A^{\al}, \al\in\Lambda\}$
will be called \textit{proper} if for each $\al\in\Lambda$ the
germ $A^{\al}$ is a proper germ ($x_{0}\in A^{\al}$).
\end{enumerate}
\end{df}

Let $(\Lambda,\ll_{\Lambda}),(\Gamma,\ll_{\Gamma})$ be directed
sets.  Let us consider $\Lambda\times \Gamma=\{(\al,\beta): \
\al\in\Lambda, \beta\in\Gamma\}$. In the product we define the
relation $\ll_{\Lambda\times\Gamma}$ in the following way:
$(\al_1,\beta_1)\ll (\al_2,\beta_2)\ \mbox{iff}\
\al_1\ll_{\Lambda}\al_2 \ \mbox{and} \
\beta_1\ll_{\Gamma}\beta_2.$

\begin{stw}
$(\Lambda\times\Gamma,\ll_{\Lambda\times\Gamma})$ is a directed
set.
\end{stw}

\begin{df}
Let  $\A=\{A^{\al},\al\in \Lambda\}, \ \B=\{B^{\beta}, \beta \in
\Gamma\}$ be systems of germs at $x_0\in K^T.$ Then we define
\begin{enumerate}
\item $\A\cap \B:=\{ A^{\al}\cap B^{\beta}, \ (\al, \beta)\in
\Lambda \times \Gamma\}.$

\item $\A\cup \B:=\{A^{\al}\cup B^{\beta}, \ (\al, \beta)\in
\Lambda \times \Gamma\}.$

\end{enumerate}
\end{df}

\begin{stw} \label{sumarodzin}
The union and the intersection of  systems of germs at $x_0$ are
systems of germs.
\end{stw}
\begin{stw}
The union and the intersection  of proper systems of germs are
proper. \mbox{}
\end{stw}

Observe that, in general,  $\A\cup \B\neq \B\cup \A$ and $\A\cap
\B\neq \B\cap \A.$ We shall introduce a relation  (between two
systems) that will allow to compare two systems in some way. In
particular we will be allowed to compare a system of germs at
$x_0$ with the point-germs system.

\begin{df}\label{porownaniek}
      Let us  consider two systems of germs at $x_0$:
     $\A=\{A^{\al},\al\in \Lambda\}$ and $\B=\{B^{\beta}, \beta \in \Gamma\}$.
     Then
     \[\A \preccurlyeq \B\Leftrightarrow \forall \beta\in \Gamma \ \exists
     \al\in \Lambda: \ A^{\al}\subseteq B^{\beta}. \]
\end{df}

It is easy to notice that the relation $\zaw$ from
Definition~\ref{porownaniek} is reflexive and transitive.

\begin{uw}
A similar relation defined for flags of finite-dimensional
algebraic varieties was used in \cite{Ty}.
\end{uw}

\begin{stw}\label{punkty}
Let  $\rp$ be the system of point-germs  at $x_0$ and
$\A=\{A^{\al},\al\in\Lambda\}$ be a proper system of germs at
$x_0$. Then $\rp \zaw \A.$
\end{stw}

The relation $\zaw$ between two systems of germs has similar
properties as the relation of inclusion of sets.

\begin{df}\label{rownosc} Let $\A,\B$ be systems of germs at $x_0$.
We define $\A\row \B\Leftrightarrow \A\zaw \B \ \mbox{and}\ \B\zaw
\A.$
\end{df}
\begin{stw}
The relation $\row$ is an equivalence relation in the collection
of systems of germs at $x_0$.
\end{stw}
\begin{df}
The equivalence class of the system $\A$  will be denoted by
$[\A]$ and called  the \textit{multigerm at $x_0$ determined by
the system $\A$.} The multigerm determined by the point-germ
system $\rp$ will be denoted by  $[\{x_0\}]$ and called  the
\textit{point-multigerm}. A trivial system determines the
\textit{empty multigerm} denoted by $[\emptyset].$
\end{df}

In the collection of multigerms  at $x_0$ the operations like
union and intersection are well-defined. Namely let $\A, \B$ be
systems of germs at $x_0.$ Then  $[\A]\cup [\B]:=[\A\cup \B],$
$[\A]\cap [\B]:=[\A\cap \B].$
 We say that $[\mathcal{A}]\subseteq
[\mathcal{B}]$ if there are $\ti{\A}\in[\A]$ and $\ti{\B}\in[\B]$
such that $\ti{\A}\zaw\ti{\B}.$

The following proposition describes a particular case of a system
of germs that determines the point-multigerm.
\begin{stw}\label{alfypunktowe}
Let $\A=\{A^{\al},\al\in\Lambda\}$ be a system of germs at $x_0$.
If for each $S\in \fin(T)$ there exists  $\al \in \Lambda$ such
that $A^{\al}=\pp(x_{0S})$ then $[\A]=[\{x_0\}_{\fin(T)}]$.
\end{stw}

The converse to the implication given in
Proposition~\ref{alfypunktowe} does not hold. We illustrate this
situation in the following example.
\begin{prz}
Let $T=\N,$ $x_0=0\in \R^{\N}$ and $\A=\{A^n,n\in\N \}$, where
$A^n=\{x\in\R^{\N}: x_1=0,\ldots,x_{n-1}=0,x_n=x_{n+1}\}_{x_0}$.
Of course the condition from the definition of the system of germs
is satisfied. We have that $\A\row \{x_0\}_{\fin(\N)},$ because
for $S\in \fin(\N)$ it is enough to take  $n=1+\max\{k\in
S\}\in\N$. Then $A^n\subset (\Pi^{\N}_S)^{-1}(x_{0S})$, but no
$A^n$ is equal to $(\Pi^{\N}_S)^{-1}(x_{0S}).$
\end{prz}

\begin{stw}\label{ostatni_system}
Let $[\A], [\B], [\mathcal{C}]$ be multigerms at $x_0$. Then
\begin{enumerate}
\item $[\A]\subseteq [\A], \ [\A]\subseteq [\A]\cup [\B], \
[\A]\cap [\B]\subseteq [\A].$ \item $[\A]\cap [\B]=[\B]\cap [\A]$
and $[\A]\cup [\B]=[\B]\cup [\A]$. \item $[\A]\cap ([\B]\cap
[\cc])=([\A]\cap [\B])\cap [\cc]$ and
     $[\A]\cup ([\B]\cup [\cc])=([\A]\cup [\B])\cup [\cc]$.
\item $[\emptyset] \subset [\A], \ [\A]\cap [\A]= [\A]$ and
$[\A]\cup [\A]= [\A]$. \item $[\A]\cap [\emptyset] =[\emptyset], \
[\A]\cup[\emptyset] =[\A].$ \item for the proper multigerm $[\A]$:
$[\A]\cap [\{x_0\}]= [\{x_0\}],$  $[\A]\cup [\{x_0\}]= [\A]$
    and $[\{x_0\}]\subseteq [\A].$
\end{enumerate}
\end{stw}

\section{Theorems of zeros}\label{roz5}

\begin{df} Let $I$ be an ideal of the ring $\On$.
The {\it zero-system of an ideal I} is defined to be the system of
germs at $x_0:$ $ \mathcal{Z}(I)=\{Z(\al),\al\in \fin(I)\},$ where
$Z(\al)$ is the zero set-germ of  $\al.$ We consider $\fin(I)$
with the inclusion relation. An equivalence class of $\cZ(I)$
(with respect to the relation $\row$) is called the
\textit{zero-multigerm of the ideal $I$} and it is denoted by
$[\cZ(I)]$.
\end{df}

\begin{stw} \label{zera} Let $I,J\lhd$  $\On$.
 If  $J\subseteq I$ then $[\cZ(I)]\subseteq [\cZ(J)].$
\end{stw}

Every proper ideal of the ring $\On$ is contained in
$m_{x_0}=(x_t-x_{0t})_{t\in T}$. Therefore:
\begin{wn}
  For a proper ideal $I$ of $\On$:
  $[\cZ(m_{x_0})]\subseteq [\cZ(I)].$
\end{wn}

\begin{stw}\label{zeramaksymalne}
  $[\cZ(m_{x_0})]=[\{x_0\}].$
\end{stw}
\begin{proof} The system $\cZ(m_{x_0})$ is proper
and from  Corollary~\ref{ostatni_system}: $[\{x_0\}]\subseteq
[\cZ(m_{x_0})].$ Now let $S\in \fin(T)$ and we take
$\al=\{x_t-x_{0t}\}_{t\in S}\in \fin(m_{x_0}).$ Then
$Z(\al)=\pp(x_{0S}).$ Therefore $[\cZ(m_{x_0})]=[\{x_0\}]$.
\end{proof}

\begin{wn}
If $I$ is a proper ideal of $\On$ then $[\cZ(I)]\neq[\emptyset].$
\end{wn}
\begin{stw}\label{sumaidealow}
Let $I,J\lhd$ $\On$. Then $[\cZ(I+J)]=[\cZ(I)]\cap [\cZ(J)].$
\end{stw}
\begin{proof} From  Proposition~\ref{zera} we have that:
$[\cZ(I+J)]\subseteq [\cZ(I)]$ and $[\cZ(I+J)]\subseteq [\cZ(J)].$
Hence $[\cZ(I+J)]\subseteq [\cZ(I)]\cap [\cZ(J)].$

Now let $\al\in \fin(I+J)$, $\al=\{\f_1,\ldots,\f_k\}$, where
$\f_i=\psi_{i,1}+\psi_{i,2}$ and  $\psi_{i,1}\in I$,
$\psi_{i,2}\in J$. Then let $\beta=(\beta_1,\beta_2)\in
\fin(I)\times \fin(J)$ be such that
$\beta_1=\{\psi_{1,1},\ldots,\psi_{k,1}\}$,
$\beta_2=\{\psi_{1,2},\ldots,\psi_{k,2}\}$. Then $Z(\beta_1)\cap
Z(\beta_2)\subseteq Z(\al).$ Hence $\cZ(I)\cap \cZ(J)\zaw
\cZ(I+J)$.
\end{proof}

\begin{stw} \label{zerarad}
Let $I$ be an ideal of $\Op_{x_0}(\R^T)$. Then
$[\cZ(I)]=[\cZ(\sqrt[\R]{I})].$
\end{stw}
\begin{proof} Because $I\subseteq \sqrt[\R]{I}$, from
Proposition~\ref{zera} we get:  $[\cZ(\sqrt[\R]{I})]\subseteq
[\cZ(I)].$ To prove  that $\cZ(I)\subseteq \cZ(\sqrt[\R]{I})$ we
must show that for any $\beta \in \fin(\sqrt[\R]{I})$ there is
$\ani $, such that $Z(\al)\subseteq Z(\beta).$ Let $\beta \in
\fin(\sqrt[\R]{I})$, hence $\beta =\{\f _1,\ldots ,\f_p\} $ and
for $i=1,\ldots, p$: $\f_i\in \sqrt[\R]{I}$. Then for each
$i=1,\ldots ,p $ there are $m_i\in \N$, $k_i\in \N \cup \{0\}$ and
$\psi _{i,1}, \ldots \psi _{i,k_i}\in \On$, such that: $\f
_i^{2m_i}+\sum \limits _{j=1}^{k_i} \psi_{i,j}^2\in I$. Let
$\al=\{\f _1^{2m_1}+\sum \limits _{j=1}^{k_1} \psi_{1,j}^2, \ldots
, \f _p^{2m_p}+\sum \limits _{j=1}^{k_p} \psi_{p,j}^2\}.$ Then
$Z(\al)\subseteq Z(\beta).$
\end{proof}

\begin{stw} \label{zerarad1}
Let $I$ be an ideal of $\Op_{x_0}(\C^T)$. Then
$[\cZ(I)]=[\cZ(\sqrt{I})].$
\end{stw}

Let $[\A]$ be a multigerm at $x_0.$ Then we define the set
$J([\A]):=\{ \f\in \On: [\A]\subseteq [\cZ(\f)]\},$
  where $[\cZ(\f)]$ is the multigerm generated by the zero-system of the ideal $(\f).$

\begin{stw}
  Let $[\A]$ be a multigerm at $x_0$
  Then $J([\A])$
  is an ideal of the ring $\On$.
\end{stw}
\begin{proof} Let $\f \in J([\A])$ and $\psi \in J([\A])$. Then
$[\A]\subseteq [\cZ(\f)]$ and $[\A]\subseteq [\cZ(\psi)]$. Hence
$[\A]\subseteq [\cZ(\f)]\cap
[\cZ(\psi)]=[\cZ((\f)+(\psi))]\subseteq [\cZ((\f+\psi))],$ (from
Propositions~\ref{sumaidealow} and \ref{zera}).

Now let $\f\in J([\A])$, $\psi\in \On$. As $(\f\cdot\psi)\subseteq
(\f)$, then from the Proposition~\ref{zera}: $[\A]\subseteq
[\cZ(\f)]\subseteq [\cZ(\f\cdot\psi)].$ Thus $\f\cdot\psi \in
J([\A]).$
\end{proof}

The ideal $J([\A])$ has similar properties to the corresponding
ideal in the finite-dimensional case.
\begin{stw}\label{idealy}
Let $[\A]$, $[\B]$ be the multigerms  at $x_0$. Then
\begin{enumerate}
\item if $[\A]\subseteq [\B]$ then $J([\B])\subseteq J([\A])$.
\item $J([\A\cup \B])=J([\A])\cap J([\B])$.
\end{enumerate}
\end{stw}

\begin{wn}
$J([\emptyset])=\On.$
\end{wn}

Let $I\lhd\On.$ Let us consider its zero-multigerm $[\cZ(I)].$
Then in $\On$ the ideal $J([\cZ(I)])=\{\f \in \On :
[\cZ(I)]\subseteq [\cZ(\f)]\}$ is well-defined. The condition
$[\cZ(I)]\subseteq [\cZ(\f)]\}$ means that  there is $\al \in
\fin(I)$ such that $Z(\al)\subseteq Z(\f).$

Now we are ready to state the main results of this paper.
\begin{tw}\label{Rislera}
Let $I$ be an ideal of $\Op_{x_0}(\R^T).$ Then
$J([\cZ(I)])=\sqrt[\R ]{I}.$
\end{tw}
\begin{proof}
First we show that $J([\cZ(I)])\subset \sqrt[\R ]{I}.$

 Let $\f \in
J([\cZ(I)])$. Then $\cZ(I)\zaw \cZ(\f).$ Hence there is $\al \in
\fin(I)$ such that $Z(\al)\subseteq Z(\f)$. From
Corollary~\ref{3.15a} there is  $S \in \fin(T)$ such that:
$Z((M^T_S)^{-1}((\al)))\subseteq Z((M^T_S)^{-1}((\f))).$ Then
$J\left(Z((M^T_S)^{-1}((\f)))\right)\subseteq
J\left(Z((M^T_S)^{-1}((\al)))\right)$
 in
$\Oh_{x_{0S}}(\R^{S})$. Now from Risler's theorem for the ring
$\Oh_{x_{0S}}(\R^{S})$ we get that
$\sqrt[\R]{(M^T_S)^{-1}((\f))}\subseteq
\sqrt[\R]{(M^T_S)^{-1}((\al))}.$ Observe now that
$(M^T_S)^{-1}(\f)\in \sqrt[\R]{(M^T_S)^{-1}((\f))}.$ Hence
$(M^T_S)^{-1}(\f)\in \sqrt[\R]{(M^T_S)^{-1}((\al))}.$ Then
\[\f=M^T_S\left((M^T_S)^{-1}(\f)\right)\in M^T_S\left(\sqrt[\R]{(M^T_S)^{-1}((\f))}\right)
\subseteq M^T_S\left(\sqrt[\R]{(M^T_S)^{-1}((\al))}\right).\] Now
from Proposition~\ref{homom1} we obtain
\[\f\in
M^T_S\left((M^T_S)^{-1}\left(\sqrt[\R]{(\al)}\right)\right)\subseteq
\sqrt[\R]{(\al)}\subseteq \sqrt[\R]{I}.\]

If $\f\in \sqrt[\R]{I}$ then $(\f)\subseteq \sqrt[\R]{I}$. Hence
from Proposition~\ref{zera}: $[\cZ(\sqrt[\R]{I})]\subseteq
[\cZ(\f)].$ In Proposition~\ref{zerarad} we showed that
 $[\cZ(\sqrt[\R]{I})]=[\cZ(I)].$ Then also $[\cZ(I)]\subseteq [\cZ(\f)]$. Therefore
 $\f \in J([\cZ(I)])$.
\end{proof}

\begin{tw}\label{Hilberta}
Let $I$ be an ideal of $\Op_{z_0}(\C^T).$ Then
$J([\cZ(I)])=\sqrt{I}.$
\end{tw}
\begin{proof} The proof is similar to the proof for
the real case. The inclusion $J([\cZ(I)])\subseteq \sqrt{I}$
follows from the R\"{u}ckert theorem \cite{Ru,Ruc} in the
finite-dimensional case.
\end{proof}

\begin{wn} \label{lokalny}
Let $I$ be an ideal of $\Oh_{x_0}(\R^T)$. Then
$[\cZ(I)]=[\{x_0\}]\ \Leftrightarrow \ \sqrt[\R]{I}=m_{x_0}.$
\end{wn}

\begin{uw}
The condition $[\cZ(I)]=[\{x_0\}]$ in the real case describes
local observability of an infinite-dimensional dynamical system.
The ideal $I$ is generated by the germs at $x_0$ of the functions
from the observation algebra of the system that vanish at $x_0$.
See \cite{Mo,MB4} for details.
\end{uw}

\begin{wn} \label{lokalny1}
Let $I$ be an ideal of $\Oh_{z_0}(\C^T)$. Then
$[\cZ(I)]=[\{x_0\}]\ \Leftrightarrow \ \sqrt{I}=m_{z_0}.$
\end{wn}

\bibliographystyle{amsplain}

\begin{thebibliography}{}

\bibitem{ABR} C. Andradas, L. Br\"{o}cker, and J.M. Ruiz,
    {\em Constructible Sets in Real Geometry}, Berlin Heidelberg New
    York: Springer-Verlag 1996.
\bibitem{Ban} S. Banach, {\em Th\'{e}orie des op\'erations
    lin\'{e}aires}, Warsaw 1932.
\bibitem{Ba1} Z. Bartosiewicz, {\em Local observability of
    nonlinear systems},
     Systems $\&$ Control Letters \textbf{25} (1995), 295--298.
\bibitem{Ba2} Z. Bartosiewicz, {\em Real analytic geometry and
    local observability}, Proc. Sympos. Pure Math. \textbf{64} (1998),
 Amer. Math. Soc., Providence, RI.
\bibitem{BCR} J. Bochnak, M. Coste, and M.-F. Roy,
    {\em G\'{e}om\'{e}trie alg\'{e}brique r\'{e}elle},
 Berlin Heidelberg New York: Springer-Verlag 1987.
\bibitem{Du} D.W. Dubois, {\em A Nullstellensatz for ordered
    fields}, Ark. Mat. \textbf{8} (1969), 111-114.
\bibitem{GR} R.C. Gunning,  H. Rossi, {\em Analytic Functions of
    Several Complex Variables},  Englewood Cliffs: Prentice-Hall 1965.
\bibitem{Kr} J.-L. Krivine, {\em Anneaux pr\'{e}ordonne},
    J.~Analyse Math. \textbf{12} (1964), 307-326.
\bibitem{La} S. Lang, {\em Hilbert's Nullstellensatz in
    infinite-dimensional space}, Proc. Amer. Math. Soc.
    \textbf{3} (1952), 407--410.
\bibitem{Mo} D. Mozyrska,  {\em Local observability of
    infinitely-dimensional finitely presented dynamical systems with
    output} (in Polish),  Ph.D. thesis, Technical University of
    Warsaw, Poland, 2000.
\bibitem{MB3} D. Mozyrska, Z. Bartosiewicz, {\em Families of germs
    in local observability of infinite--dimensional dynamical
    systems}, in: Proceedings of the First International Conference on
    Control and Self--Organization in Nonlinear Systems (Ed.:
    Z.~Bartosiewicz, M.~Marczak, E.~Paw{\l}uszewicz), Bia\l ystok
    (Supra\'sl), Poland, 15--18 February 2000.
\bibitem{MB4} D. Mozyrska, Z. Bartosiewicz, {\em Local observability
    of systems on $\R^{\infty}$}, in: Proceedings of Internatinal
    Conference MTNS 2000 (Ed.: M.~Fliess), Perpignan, France, June
    2000.
\bibitem{Ras} H. Rasiowa, {\em Introduction to Modern Mathematics},
    American Elsevier Publishing Co.,Inc., New York 1973.
\bibitem{Ri0} J.-J. Risler, {\em Une caract\'{e}risation des
    id\'{e}aux des vari\'{e}t\'{e}s alg\'{e}briques r\'{e}elles},
    C.R. Acad. Sci. Paris \textbf{271} (1970), 1171--1173.
\bibitem{Ri} J.-J. Risler, {\em Le th\`{e}or\'{e}me des z\`{e}ros en
    g\`{e}om\`{e}tries alg\`{e}brique   et analytique r\`{e}elles},
     Bull. Soc. Math. France \textbf{104} (1976),  113--127.
\bibitem{Ruc} W. R\"{u}ckert, {\em Zum Eliminationsproblem der
    Potenzenreihenideale}, Math. Ann. \textbf{107} (1932),
    259--181.
\bibitem{Ru} J. Ruiz,  {\em The Basic Theory of Power Series},
    Vieweg 1993.
\bibitem{Ty} A.N. Tyurin,  {\em Finite-dimensional bundles on infinite
    varieties} (in Russian), Izv. Akad. Nauk SSSR Ser. Mat. \textbf{40}(1976), 1248--1268.
\end{thebibliography}

\textbf{Authors' affiliation:}\\
 \author{Dorota Mozyrska \& Zbigniew Bartosiewicz \\
             Institute of Mathematics and Physics\\
            Bia\l ystok Technical University \\ Wiejska 45A, 15-351 Bia\l ystok,
            Poland\\ \texttt{admoz@w.tkb.pl, bartos@cksr.ac.bialystok.pl}}

          \end{document}